\def\year{2020}\relax
\documentclass[letterpaper]{article} % DO NOT CHANGE THIS
\usepackage{aaai20}  % DO NOT CHANGE THIS
\usepackage{times}  % DO NOT CHANGE THIS
\usepackage{helvet} % DO NOT CHANGE THIS
\usepackage{courier}  % DO NOT CHANGE THIS
\usepackage[hyphens]{url}  % DO NOT CHANGE THIS
\usepackage{graphicx} % DO NOT CHANGE THIS
\usepackage{graphicx}  % DO NOT CHANGE THIS
\usepackage{amsmath}
\usepackage{bbm}
\usepackage{cleveref}
\crefname{equation}{Eq.}{Eqns.}
\crefname{figure}{Fig.}{Figs.}
\crefname{table}{Table}{Tables}
\title{Solving Inverse Problems in Steady-State Navier-Stokes Equations using Deep Neural Networks} % NN, steady-state Navier-Stokes Equations, Inverse Problem

\author{Tiffany Fan,\textsuperscript{\rm 1}%\thanks{}  

Kailai Xu,\textsuperscript{\rm 1}
Jay Pathak,\textsuperscript{\rm 2}  
Eric Darve\textsuperscript{\rm 1, \rm 3} \\ % All authors must be in the same font size and format. Use \Large and \textbf to achieve this result when breaking a line
\textsuperscript{\rm 1}Institute for Computational and Mathematical Engineering, Stanford University, Stanford, CA 94305, USA; \\ %If you have multiple authors and multiple affiliations
\{tiffan, kailaix, darve\}@stanford.edu  \\% email address must be in roman text type, not monospace or sans serif
\textsuperscript{\rm 2} Ansys Inc., San Jose, CA 95134, USA; jay.pathak@ansys.com \\
\textsuperscript{\rm 3} Mechanical Engineering, Stanford University, Stanford, CA 94305, USA
}

\begin{document}

\maketitle

\begin{abstract}
Inverse problems in fluid dynamics are ubiquitous in science and engineering, with applications ranging from electronic cooling system design to ocean modeling. We propose a general and robust approach for solving inverse problems in the steady-state Navier-Stokes equations by combining deep neural networks and numerical partial differential equation (PDE) schemes. Our approach expresses numerical simulation as a computational graph with differentiable operators. We then solve inverse problems by constrained optimization, using gradients calculated from the computational graph with reverse-mode automatic differentiation. This technique enables us to model unknown physical properties using deep neural networks and embed them into the PDE model. We demonstrate the effectiveness of our method by computing spatially-varying viscosity and conductivity fields with deep neural networks (DNNs) and training the DNNs using partial observations of velocity fields. We show that the DNNs are capable of modeling complex spatially-varying physical fields with sparse and noisy data. Our implementation leverages the open access ADCME, a library for solving inverse modeling problems in scientific computing using automatic differentiation.
\end{abstract}

\section{Introduction} \label{Sec:intro}

% Figure: computational graph 

% Algorithm or Figure:  fluid solver

% Figure:  to compare NN vs. Variables

Fluid dynamics is fundamental for a wide variety of applications in aeronautics, geoscience, meteorology and mechanical engineering, such as chip design \cite{fedorov2000three}, earth exploration \cite{li2020coupled}, and weather forecasting \cite{zajaczkowski2011preliminary}. However, quantifying fluid properties in the governing equations, which are essential for predictive modeling, remains a challenging problem: the computation can be expensive %or impossible to determine from direct experiments 
and often leads to underdetermined or ill-posed systems \cite{cotter2009bayesian}. This challenge leads us to leverage indirect data, which are not direct observations of fluid properties but information related to the fluid properties via the governing equations. To utilize the indirect data, we need to consider the relationship between governing equations, data, and fluid properties as a whole.

Due to the high computational cost of high-fidelity numerical simulations \cite{freund2019dpm}, researchers have made significant efforts in leveraging machine learning to assist fluid dynamics simulations in the past decade. Based on simulation solutions of Reynolds averaged Navier-Stokes (RANS) models, Ling and Templeton trained classifiers to infer the error-prone regions of the domain \cite{ling2015evaluation}. Other researchers utilized deep neural networks to learn complex physical relationships to predict the simulation outcomes of fluid dynamics systems, where the neural networks are used either as approximations for implicit functions to provide end-to-end predictions for simulation outcomes \cite{sirignano2018dgm} or as augmentations to (simplified or partially known) physical laws to predict correction terms for the simulation outcomes \cite{wang2017physics,wu2017priori,freund2019dpm}. Holland, Baeder, and Duraisamy integrated neural networks into the Field Inversion and Machine Learning (FIML) approach to improve RANS predictions for airfoils by learning the spatially-varying discrepancy corrections  \cite{holland2019field}.

The governing equations considered in this paper are the steady-state Navier-Stokes equations for incompressible flow
\begin{equation}\label{eq:mass_and_momentum}
\begin{aligned}
    (\mathbf{u} \cdot \nabla) \mathbf{u} =&
    -\frac{1}{\rho} \nabla p + \nabla\cdot (\nu \nabla \mathbf{u}) + \mathbf{g}\\
    \nabla \cdot \mathbf{u} =& \  0
\end{aligned}
\end{equation}
%with $\nabla \cdot \mathbf{u} = 0$, and 
where $\rho$ is the fluid density, $\mathbf{u}$ is the vector of flow velocity, $p$ is the pressure, $\nu$ is the kinematic viscosity field, and $\mathbf{g}$ is the vector of body accelerations. 
The system is highly nonlinear and thus the numerical simulation requires Newton's iterations. 

One example of an inverse problem is to estimate a spatially-varying viscosity field $\nu(\mathbf{x})$ from partially observed velocity data $\mathbf{u}$. Here the observations are considered ``indirect'' because $\nu(\mathbf{x})$ is not directly measured. Due to the limited number of observations, the inverse problem may be ill-posed, i.e., there may be multiple functions $\nu(\mathbf{x})$ that produce the same set of observations. Here, we propose using a deep neural network (DNN) \cite{goodfellow2016deep} to approximate the fluid properties of interest, such as $\nu(\mathbf{x})$, as a regularizer. The inputs of the DNN are geometrical coordinates and the outputs are the values of $\nu(\mathbf{x})$ at the corresponding locations. 

The main contribution of this work is to propose a general approach to couple DNNs with iterative partial differential equation (PDE) solvers for the steady-state Navier-Stokes equations. The basic idea is to express both the DNNs and the PDE solvers with a computational graph. Therefore, once we implement the forward simulation, the numerical gradients can be easily extracted from the computational graph using reverse-mode automatic differentiation \cite{margossian2019review,baydin2017automatic}. The key is to design and implement a collection of numerical simulation operators with the capability of gradient back-propagation. 

We present our new approach with its mathematical intuitions and implementation details in Sec. \ref{Sec:method}.
In Sec. \ref{Sec:numerical}, we illustrate our method with three fluid dynamics systems involving the steady-state Navier-Stokes equations and discuss the advantages of using DNNs in our method. In Sec. \ref{SubSec:numerical2} and  Sec. \ref{SubSec:numerical3}, the Navier-Stokes equations are coupled with the heat equation and the transport equation, respectively. We summarize the limitations and ongoing work of our approach in Sec. \ref{Sec:discussion}.

\section{Methodology} \label{Sec:method}
%Inverse problems aim at exploring the causal factors of given observations. %under the assumption of certain physical rules. 

%The goal of inverse learning is to predict unknown model parameters that are consistent with both the observed data and the physical laws described by the governing equations. In science and engineering applications, the model parameters represent physical properties such as material properties and initial conditions. These physical properties of interest are either constants or high dimensional and complex functions over the domain. The latter case is more challenging and thus is the focus of this work. 

%We propose the method of physics constrained machine learning to construct flexible approximations to complex model parameters without violating the physical constraints. 

\subsection{Constrained Optimization for Inverse Problems}

In this work, 
we formulate an inverse problem in the context of constrained optimization. For a system of PDEs with solution $u$ where the physical properties of interest are represented as functions of variables $\theta$, we note that $u$ is indirectly a function of $\theta$ \cite{xu2020physics}. Thus, we consider the constrained optimization problem
% \begin{gather*}
%     \min_\theta \; L(u) \; \text{($u$ is indirectly a function of $\theta$)} \\ 
%     \text{s.t.} \; F(u, \theta) = 0
% \end{gather*} 
\begin{equation}\label{eq:constrained_opt}
\begin{aligned}
    \min_\theta  & \quad  L(u) \\
    \text{s.t.}  & \quad F(u, \theta) = 0
\end{aligned}
\end{equation}

\noindent where (i) the variables $\theta$ are model parameters or neural network weights and biases; (ii) the objective function $L$ for minimization is a loss function that measures the discrepancy between the observations and the predictions of $u$ generated by forward computation; (iii) the constraints $F$ are the governing equations derived from physical laws. By solving \cref{eq:constrained_opt}, we obtain the optimal values for the variables, which provide the optimal approximation to the physical properties of interest. 

In our method, the constrained optimization problem is converted to an unconstrained problem. Specifically, we first solve the governing equations $F(u, \theta)=0$ to obtain $u = G(\theta)$. Then, we solve the unconstrained optimization problem
\begin{equation}\label{eq:unconstrained_opt}
    \min_\theta \; L(G(\theta))
\end{equation}

\subsection{Deep Neural Networks for Inverse Problems}

% Figure: DNN + Inverse Problem 

When the physical properties of interest are complex and spatially-varying, the resulting optimization problems have infinite-dimensional feasible spaces. Solving such optimization problems with traditional basis functions, such as piecewise linear basis functions and radial basis functions, is difficult due to the curse of dimensionalities and the uneven distribution of observations \cite{huang2020learning}. Thus, we use deep neural networks (DNNs) to approximate the physical properties to ensure the flexibility of the approximation. The inputs of the DNNs are spatial coordinates of a given point in the domain, and the outputs are the predicted physical properties at that point. DNNs are high dimensional nonlinear functions of the inputs and have demonstrated abilities to approximate complex unknown functions. Besides, with certain choices of activation functions, such as tanh and sigmoid functions, the DNNs are continuous functions of the inputs. The inherent regularity effect of the DNNs results in more accurate approximations to the true physical properties in many applications \cite{ulyanov2020deep}. 

\begin{figure}[htbp]
    \centering
    \includegraphics[width=\linewidth]{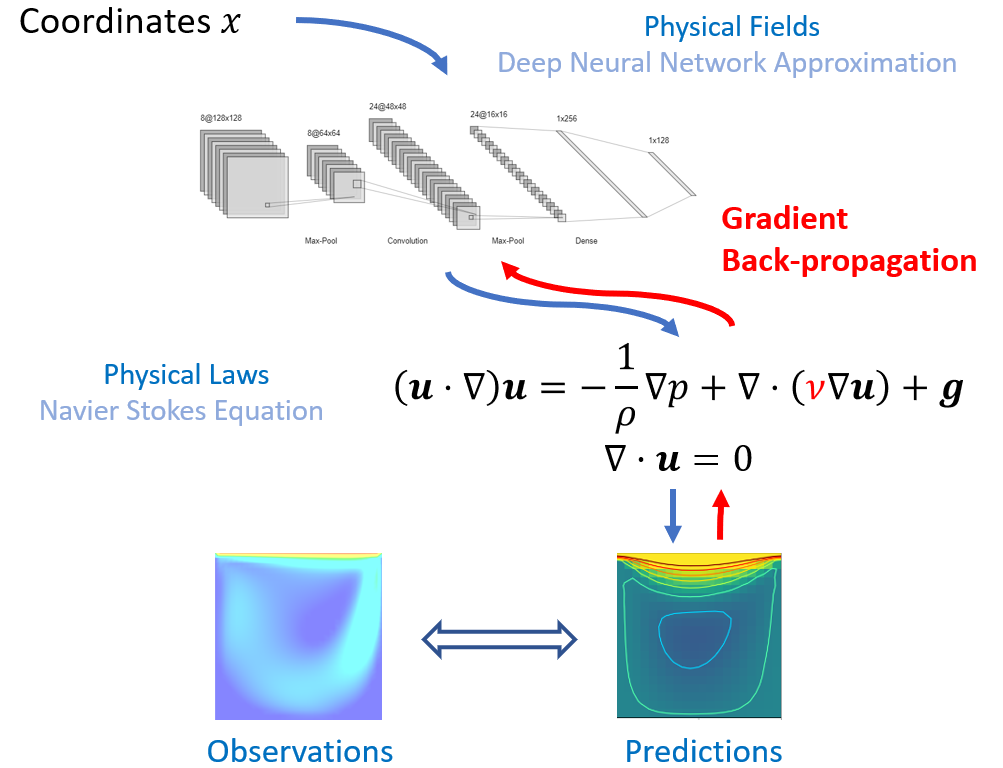}
    \caption{Schematic diagram of using DNNs to approximate spatially-varying physical properties. Note that only $\nu(\mathbf{x})$ is approximated by a DNN. The predictions of $u$ are defined on the discretized grid points and computed with the PDE solvers. The DNNs are coupled with PDE solvers and the gradients with respect to the loss function are back-propagated through both the DNNs and the PDE solvers .}
    \label{fig:advertisement}
\end{figure}

\subsection{Expressing Numerical Simulation using a Computational Graph}

\begin{figure}[htbp]
% \vspace{-0.03\linewidth}
    \centering
    \includegraphics[width=\linewidth]{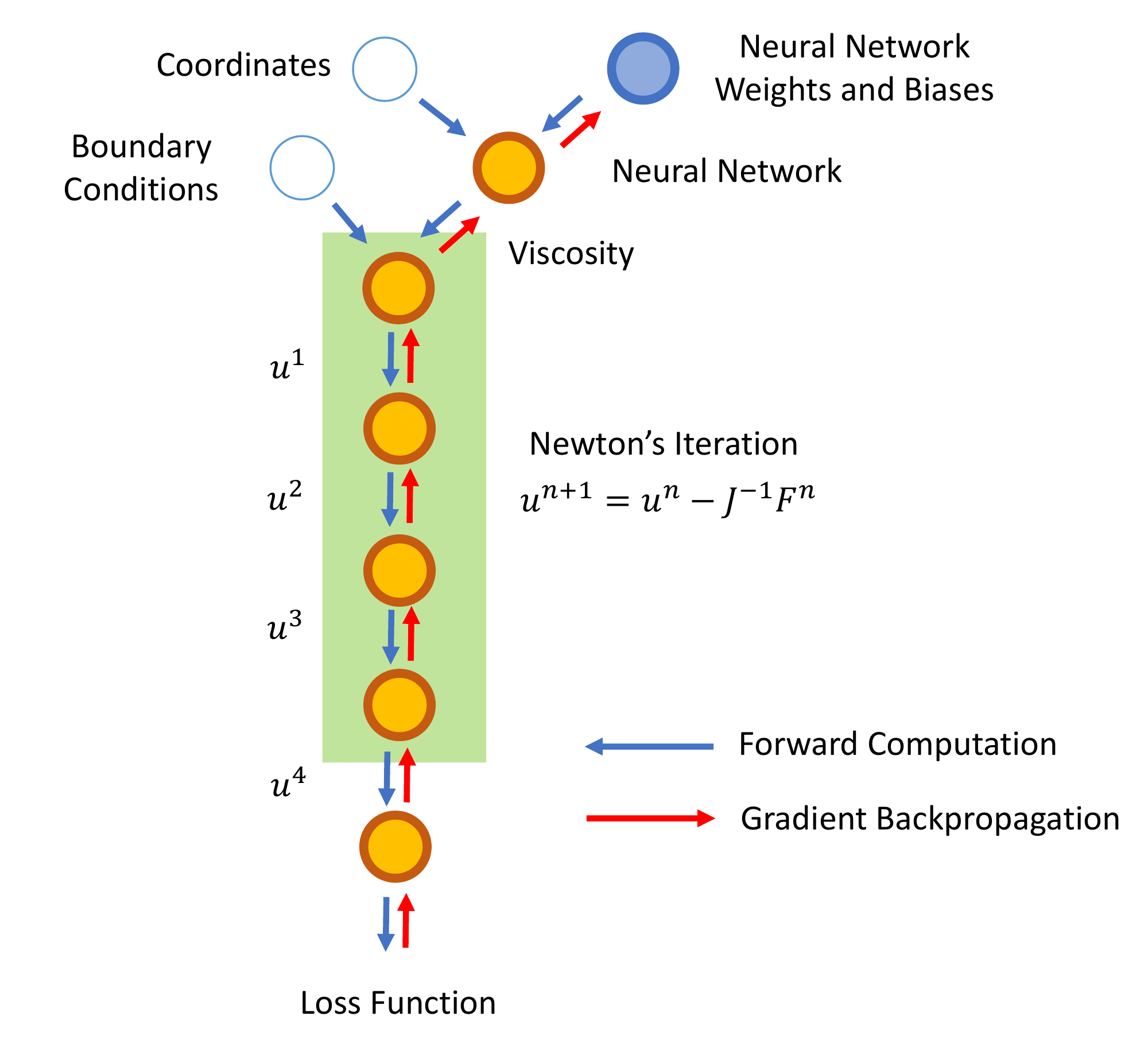}
    \caption{Expressing numerical simulation as a computational graph. The orange nodes denote numerical operators. The solid blue node denotes variables that are updated during the process of optimization. The empty blue nodes denote parameters that are fixed in this process.}
    \label{fig:computational_graph}
\end{figure}

To solve the inverse problem formulated in \cref{eq:unconstrained_opt}, we use a gradient-based optimization algorithm. The critical step is to compute the gradients $\nabla_\theta L(G(\theta))$. We represent the numerical simulaton with a computational graph with differentiable operators, as shown in  \cref{fig:computational_graph}. In each iteration of the optimization algorithm, we first perform forward computation to solve the governing equations based on the current DNN weights and biases. Next, we evaluate the loss function by comparing the computed physical quantities with the observed data. Then, we compute the gradients using reverse-mode automatic differentiation. Finally, we update the DNN weights and biases according to the optimization algorithm, using the numerical gradients. In the numerical examples, we use the L-BFGS-B optimization algorithm \cite{liu1989limited}, which performs a line search in the direction of gradient descent in every iteration.

% Figure: Example

\subsection{Physics Constrained Learning for the Nonlinear Fluid Solver}

\Cref{eq:mass_and_momentum} describes the motion of Newtonian fluids. The corresponding weak form of the steady-state Navier-Stokes equations is given by \cite{rannacher2000finite}
\begin{align}
   &\left(u \frac{\partial u}{\partial x}, \delta u'\right)+ \left(v \frac{\partial u}{\partial y} , \delta u'\right)\notag \\ 
   &  \qquad =  \frac{1}{\rho} \left(p, \frac{\partial \delta u'}{\partial x}\right)-\left(\nabla u, \nu\nabla \delta u'\right)+\left(f, \delta u'\right)\label{equ:weak1} \\ 
    &\left(u \frac{\partial v}{\partial x}, \delta v'\right) + \left(v \frac{\partial v}{\partial y}, \delta v'\right)\notag \\ 
&   \qquad   = \frac{1}{\rho}\left(p, \frac{\partial \delta v'}{\partial y}\right) -\left(\nabla v, \nu \nabla\delta v'\right)+ \left( g, \delta v'\right)\label{equ:weak2} \\ 
   & \left(\frac{\partial u}{\partial x}, \delta p'\right) + \left(\frac{\partial v}{\partial y}, \delta p' \right)  = 0\label{equ:weak3}
\end{align}
where $\mathbf{u} = (u, v)$ and $p$ are the trial functions, $\delta u'$, $\delta v'$ and $\delta p'$ are the corresponding test functions; $\nu$ is the viscosity field. Note that the system  (\Cref{equ:weak1,equ:weak2,equ:weak3}) is highly nonlinear. To solve the nonlinear system, we use Taylor's expansion to linearize the equation, before applying the Newton's iterative method \cite{beam1988newton}.

The idea is to construct a computational graph and express all the computation using differentiable operators. For example, the term $\left(\nabla u, \nu\nabla \delta u'\right)$ corresponds to a sparse block in the Jacobian matrix. We need an operator that consumes $\nu$ and outputs the sparse block. The operator should also be able to back-propagate downstream gradients, i.e., to compute $\frac{\partial L(a(\nu))}{\partial \nu}$ given $\frac{\partial L(a)}{\partial a}$, where $L$ is the scalar loss function and $a$ represents the entries in the sparse block. We refer readers to \cite{xu2020physics} on how to derive and implement such operators using physics constrained learning (PCL). 

The numerical solver for the Navier-Stokes equation is iterative. However, we found that the solver typically converges very fast (e.g., within 5 iterations for the numerical experiments in Sec. \ref{Sec:numerical}) to a very small residual. Therefore, in the gradient back-propagation, we can differentiate through each iteration in the forward computation, as is shown in \Cref{fig:advertisement}.

\section{Numerical Experiments} \label{Sec:numerical}

In this section, we apply our method to three fluid dynamics systems involving the steady-state Navier-Stokes equations. In all three examples, the DNNs share the same architecture: 3 fully connected layers, 20 neurons per layer, and \texttt{tanh} activation functions. Our implementation leverages the open access ADCME library.  

\subsection{Learning Spatially-Varying Viscosity in Steady-State Navier-Stokes Equations} \label{SubSec:numerical1}

We evaluate our method with the classic lid-driven cavity flow problem. The governing equation is given by \Cref{eq:mass_and_momentum}  with the a spatially-varying viscosity field $\nu(x, y)$.
% \begin{equation*}
%     \nu(x, y)=1+6 x^{2}+\frac{x}{1+2 y^{2}}
% \end{equation*}
We approximate $\nu(x, y)$ with a DNN, denoted $\nu_\theta(x, y)$, where $\theta$ is the weights and biases of the DNN. 

\begin{figure}[bth] 

%\hspace{.01\linewidth}

\begin{minipage}[m]{0.32 \linewidth}
\centerline{\small{reference}}
\end{minipage}
\begin{minipage}[m]{0.32\linewidth}
\centerline{\small{DNN estimation}}
\end{minipage}
\begin{minipage}[m]{0.32\linewidth}
\centerline{~~\small{pointwise estimation}}
\end{minipage}\\
% \vspace{-0.05\linewidth}
\begin{minipage}[m]{0.32\linewidth}
\centerline{\includegraphics[width=1.1\linewidth]{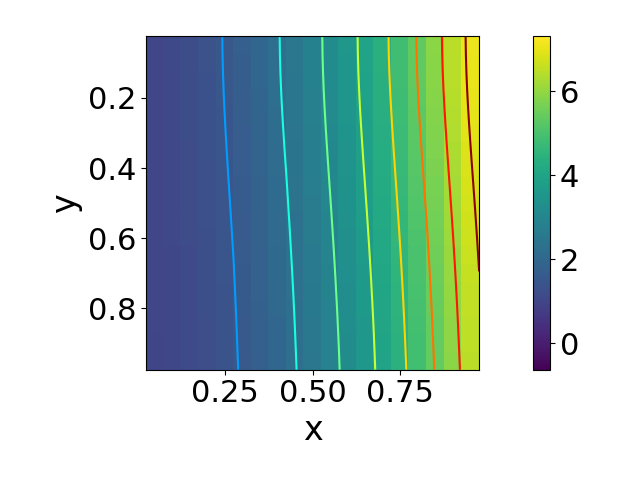}}
\end{minipage}
\begin{minipage}[m]{0.32\linewidth}
\centerline{\includegraphics[width=1.1\linewidth]{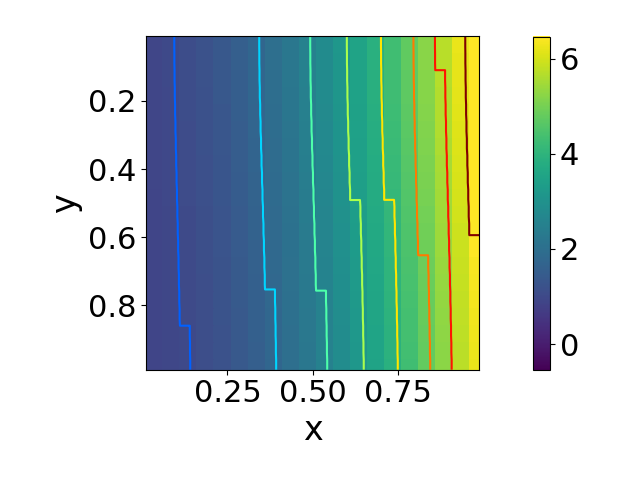}}
\end{minipage}
\begin{minipage}[m]{0.32\linewidth}
\centerline{\includegraphics[width=1.1\linewidth]{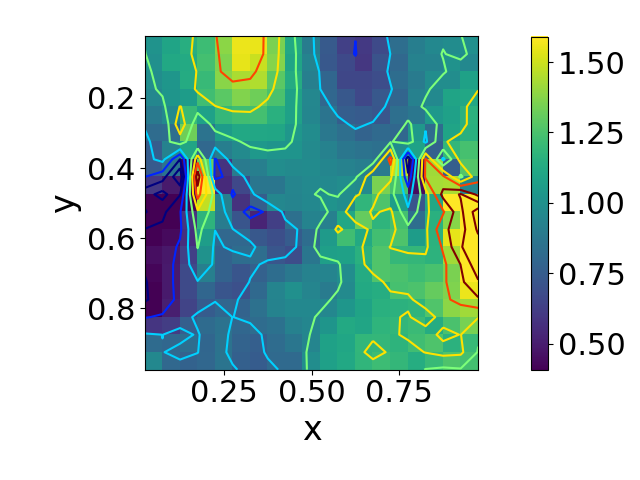}}
\end{minipage}\\
\caption{A comparison of estimations by the DNN and by pointwise values. The DNN provides a continuous and more accurate approximation to the reference viscosity function. }\label{Fig:nn_vs_pointwise}
\end{figure}

In this example, the observations are $\mathbf{u}$ at grid points, and the pressure is unknown. The observations are simulated on a grid of size $21 \times 21$, with constant density $\rho=1 $,  velocity $u=\mathbbm{1}_{\{ y=0\}}$ and $v=0$, and viscosity field
\begin{equation*}
    \nu(x, y)=1+6 x^{2}+\frac{x}{1+2 y^{2}}
\end{equation*} 
We use the velocity data to train the DNN $\nu_\theta (x,y)$. The pressure is assumed to be unknown. The estimated viscosity field is shown in \Cref{Fig:nn_vs_pointwise}. 

% Figure: Pixel vs Neural Network

\begin{figure}[htbp] 

%\hspace{.01\linewidth}

\begin{minipage}[m]{0.33 \linewidth}
\centerline{~~~\small{reference}}
\end{minipage}
\begin{minipage}[m]{0.32\linewidth}
\centerline{\small{prediction}}
\end{minipage}
\begin{minipage}[m]{0.32\linewidth}
\centerline{\small{difference}~~}
\end{minipage}
\\
% velocity x
\begin{minipage}[m]{0.01\linewidth}
\centerline{\small{$u$}}
\end{minipage}
\begin{minipage}[m]{0.32\linewidth}
\centerline{~~\includegraphics[width=1.05\linewidth]{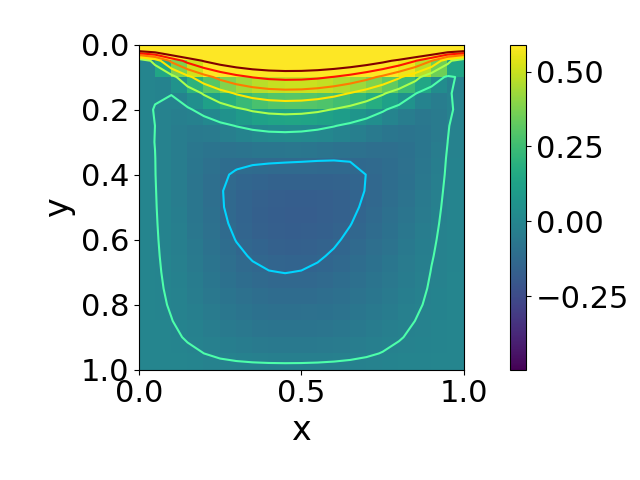}}
\end{minipage}
\begin{minipage}[m]{0.32\linewidth}
\centerline{\includegraphics[width=1.05\linewidth]{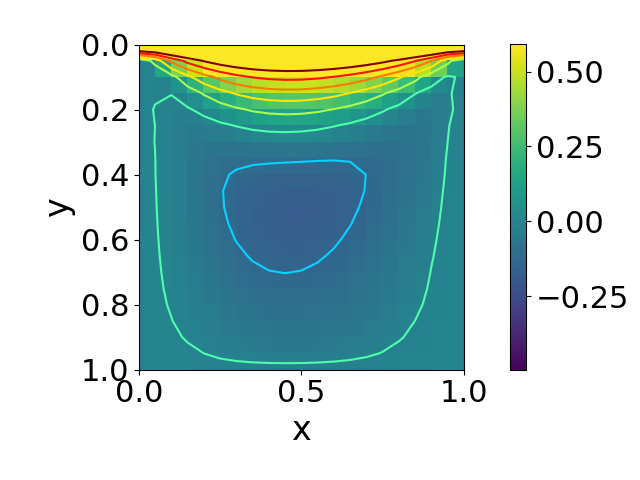}}
\end{minipage}
\begin{minipage}[m]{0.32\linewidth}
\centerline{~~~\includegraphics[width=1.15\linewidth]{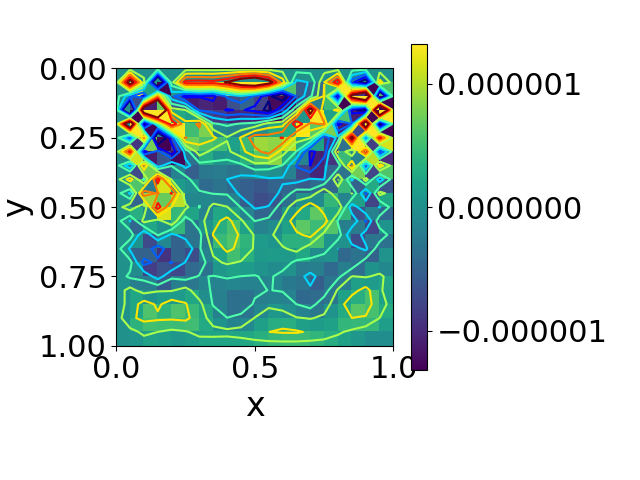}}
\end{minipage}\\
% velocity y
\begin{minipage}[m]{0.01\linewidth}
\centerline{\small{$v$}}
\end{minipage}
\begin{minipage}[m]{0.32\linewidth}
\centerline{~~\includegraphics[width=1.05\linewidth]{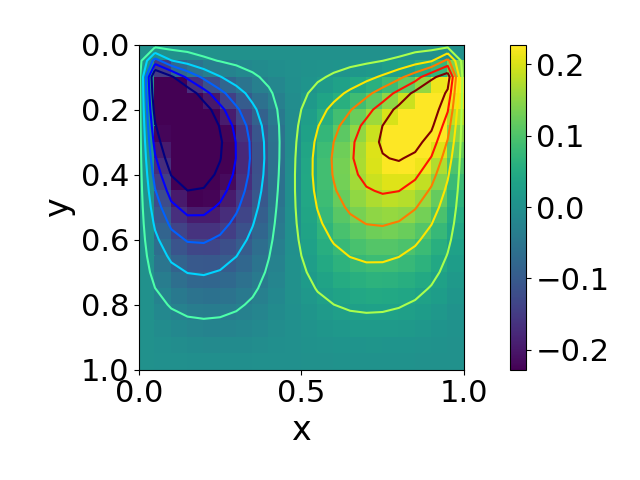}}
\end{minipage}
\begin{minipage}[m]{0.32\linewidth}
\centerline{\includegraphics[width=1.05\linewidth]{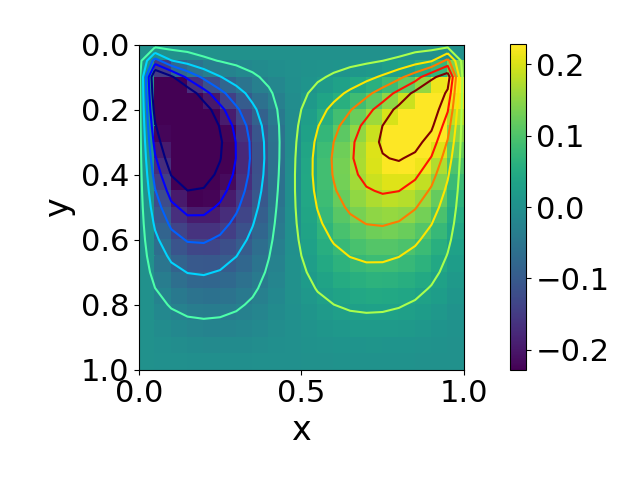}}
\end{minipage}
\begin{minipage}[m]{0.32\linewidth}
\centerline{~~~\includegraphics[width=1.15\linewidth]{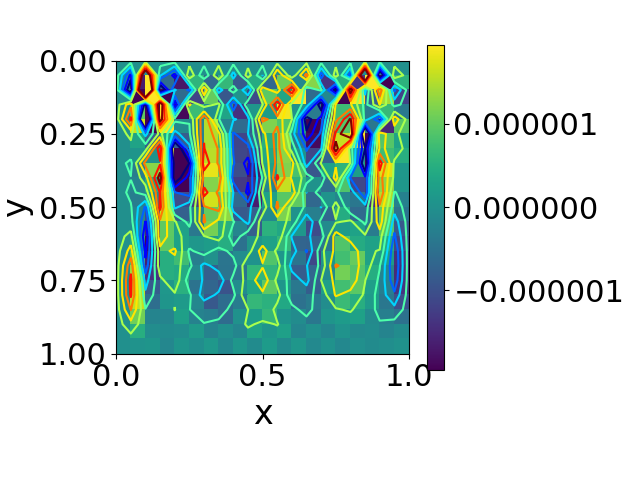}}
\end{minipage}\\
% pressure
\begin{minipage}[m]{0.01\linewidth}
\centerline{\small{$p$}}
\end{minipage}
\begin{minipage}[m]{0.32\linewidth}
\centerline{~~\includegraphics[width=1.05\linewidth]{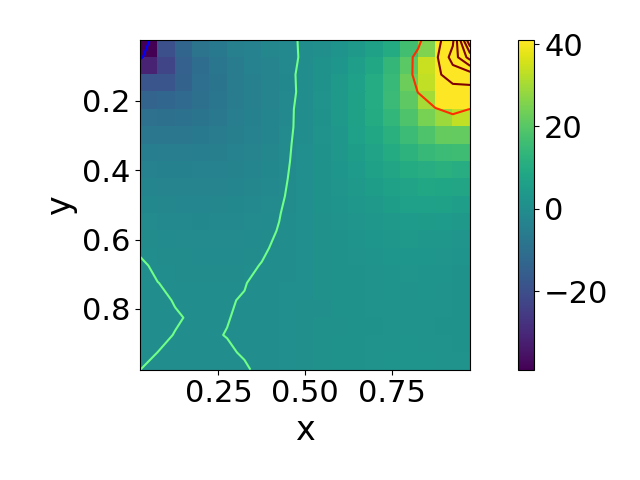}}
\end{minipage}
\begin{minipage}[m]{0.32\linewidth}
\centerline{\includegraphics[width=1.05\linewidth]{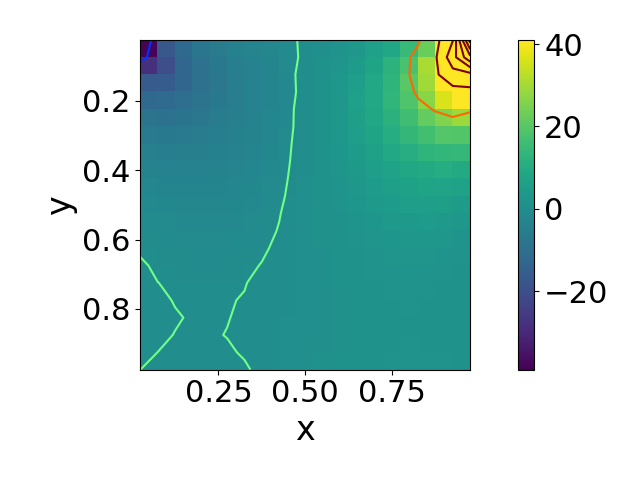}}
\end{minipage}
\begin{minipage}[m]{0.32\linewidth}
\centerline{\includegraphics[width=1.05\linewidth]{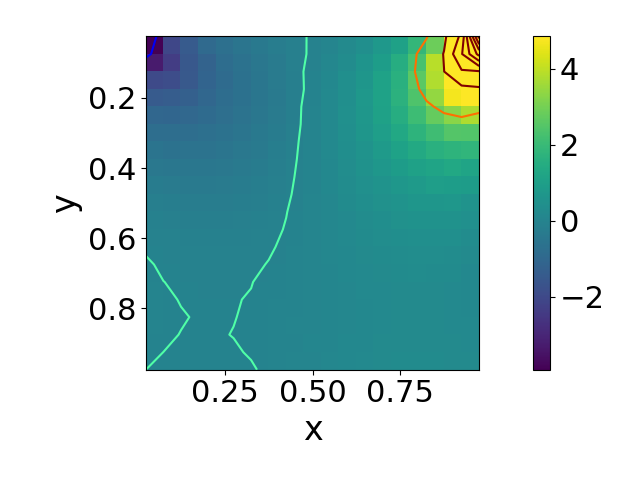}~~}
 \end{minipage} %\vspace{-0.08\linewidth}
\\
\caption{The reference velocity and pressure, the predictions from the DNN, and the corresponding error. We note that the loss function only contains the velocity data. However, benefiting from the physical constraints imposed by the numerical schemes, the DNN provides an accurate prediction of the pressure as well. \label{Fig:eg1}}
\end{figure}

\noindent Then, we plug the estimated viscosity field into \Cref{eq:mass_and_momentum} and solve for $u, v$, and $p$. The results are shown in \Cref{Fig:eg1}. We observe that the predictions are very close to the reference. 

% \vspace{-0.05\linewidth}
We also compare the DNN results with those of pointwise estimation, where we optimize the values of viscosity $\nu(x,y)$ at each grid point instead of using a function approximator (e.g., DNN). In \Cref{Fig:nn_vs_pointwise}, We observe that the DNN produces a smooth profile of the viscosity field, with a relative mean square error of $1.26$\%. The pointwise estimation is far from the exact $\nu(x,y)$, with a relative mean square error of $59.14$\%, despite producing velocity predictions that are similar to the observations. 

\begin{figure}[htbp] 
\centerline{\includegraphics[width=0.65\linewidth]{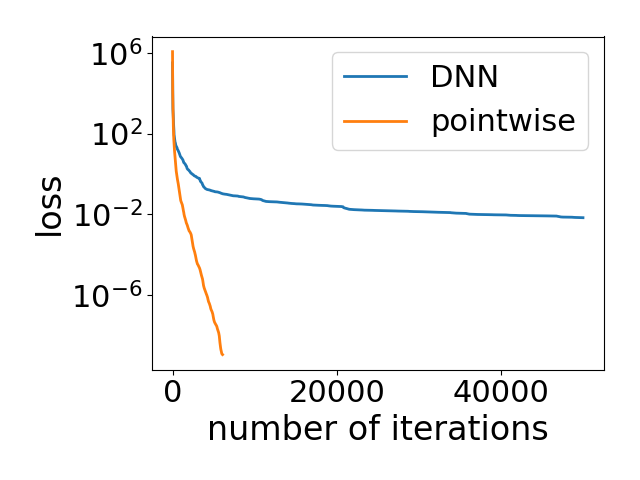}}

\caption{A comparison of the loss convergence for the DNN estimation and the pointwise estimation. }\label{fig:eg2_loss}
\end{figure}

\Cref{fig:eg2_loss} shows the convergence of the loss functions for both the DNN estimation and the pointwise estimation. The pointwise estimation achieves a smaller loss because the viscosity representation of the pointwise estimation is less constrained than that of the DNN. However, the DNN provides a better viscosity estimation due to its regularization effect.

\begin{figure}[htpb] 

% %\hspace{.01\linewidth}

\begin{minipage}[m]{0.48\linewidth}
\centerline{~~\small{DNN difference}}
\end{minipage}
\begin{minipage}[m]{0.48\linewidth}
\centerline{\small{pointwise difference}~}
\end{minipage}\\

% \vspace{-0.3cm}

\begin{minipage}[m]{0.48\linewidth}
\centerline{~~\includegraphics[width=1.05\linewidth]{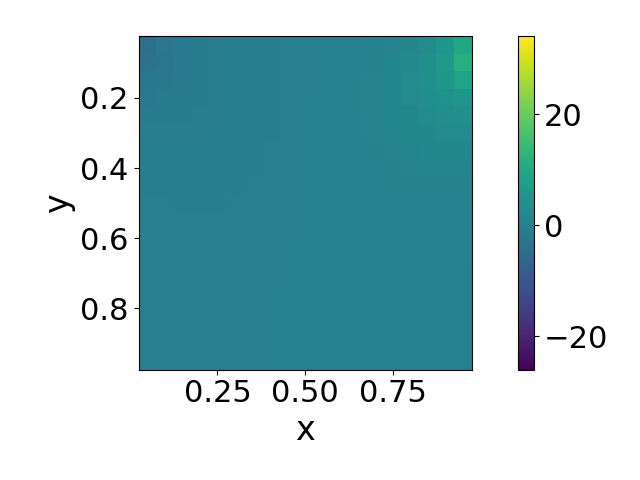}}
\end{minipage}
\begin{minipage}[m]{0.48\linewidth}
\centerline{~~\includegraphics[width=1.05\linewidth]{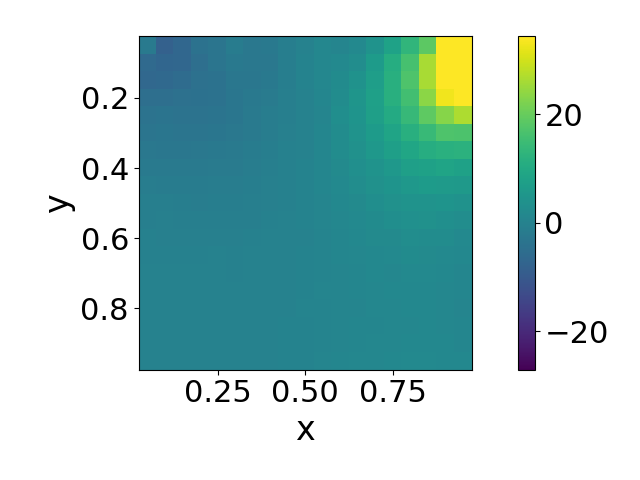}}
\end{minipage}\\

\caption{A comparison of the error in pressure predictions for the DNN estimation and the pointwise estimation. }\label{fig:eg1_pressure_error}
\end{figure}

Further evidence is shown in \cref{fig:eg1_pressure_error}, where we compare the error in the pressure predictions provided by the DNN and the pointwise estimation. The pressure profile is unobserved by both methods. Thus, for the pointwise estimation, the large error in pressure prediction and the small training loss indicate potential overfitting of the observed data. On the other hand, the DNN estimation produces an accurate estimation of the real pressure field without observing the pressure profile, due to the regularization effect of DNNs. 

\subsection{Learning Spatially-Varying Conductivity in Conjugate Heat Transfer Navier-Stokes Equations} \label{SubSec:numerical2}

In this example, we consider the coupled system with \Cref{eq:mass_and_momentum} and the energy equation (heat equation):
\begin{equation}
\label{eq:energy}
    \rho C_{p} \mathbf{u} \cdot \nabla T =
    \nabla \cdot(k \nabla T) + Q
\end{equation}
where $C_{p}$ is the specific heat capacity, $T$ is the temperature, $k$ is the conductivity, and $Q$ is the power source. The problem arises from conjugate heat transfer analysis \cite{wang2007lattice}, where the heat transfers between solid and fluid domains by exchanging thermal energy at the interfaces between them. 

\begin{figure}[htpb] 

%\hspace{.01\linewidth}

\begin{minipage}[m]{0.32 \linewidth}
\centerline{\small{reference}}
\end{minipage}
\begin{minipage}[m]{0.32\linewidth}
\centerline{\small{estimation}~~~~}
\end{minipage}
\begin{minipage}[m]{0.32\linewidth}
\centerline{\small{difference}~~~}
\end{minipage}\\

\begin{minipage}[m]{0.32\linewidth}
\centerline{\includegraphics[width=1.05\linewidth]{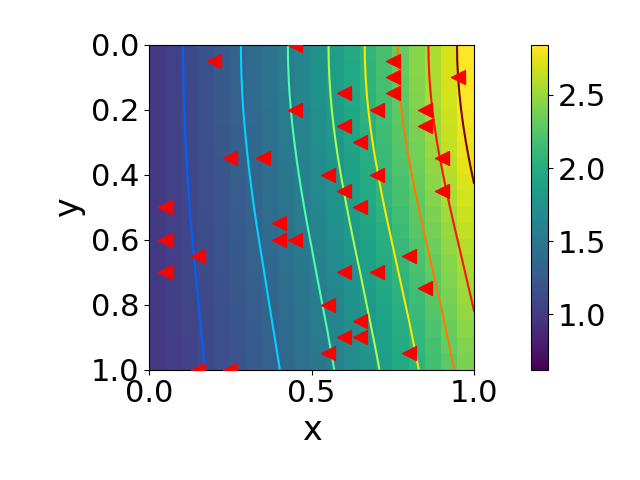}}
\end{minipage}
\begin{minipage}[m]{0.32\linewidth}
\centerline{\includegraphics[width=1.05\linewidth]{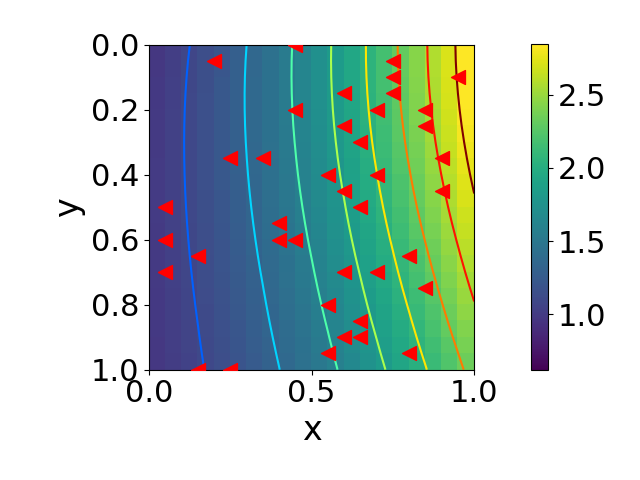}~~~~}
\end{minipage}
\begin{minipage}[m]{0.32\linewidth}
\centerline{\includegraphics[width=1.05\linewidth]{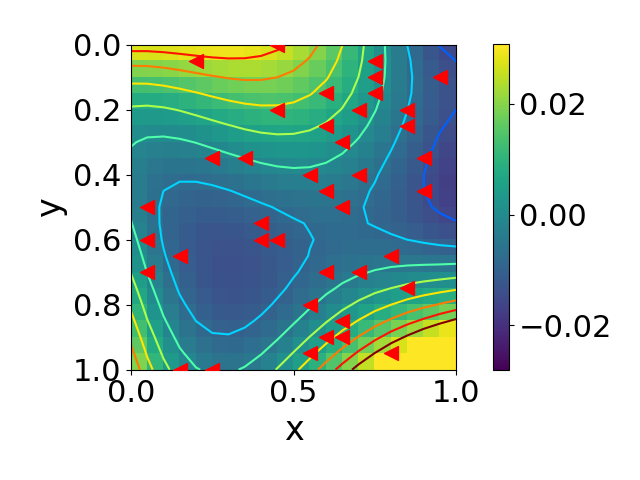}}
\end{minipage}\\
\caption{The reference conductivity, the estimated conductivity by the DNN, and the estimation error after 100 optimization steps. The 40 randomly sampled grid points, where velocity and temperature data are observed, are labeled with red triangles. }\label{fig:eg2_boussinesq}
\end{figure}

We simulate the velocity, pressure and temperature data via forward computation with $\rho=1$, $C_p=1,$ and conductivity field
\begin{equation*}
    k(x,y) = 1 + x^2 + \frac{x }{1+y^2}
\end{equation*}
The observations are the velocity and temperature data at 40 randomly sampled locations from the $21 \times 21$ grid. The pressure is assumed to be unknown. Fig.\ref{fig:eg2_boussinesq} shows that the DNN produces an accurate approximation for conductivity from the limited data.

\begin{figure}[htpb] 

%\hspace{.01\linewidth}

\begin{minipage}[m]{0.12 \linewidth}
\centerline{\small{$\epsilon$}}
\end{minipage}
\begin{minipage}[m]{0.43\linewidth}
\centerline{~~\small{estimation}}
\end{minipage}
\begin{minipage}[m]{0.43\linewidth}
\centerline{\small{difference}~}
\end{minipage}\\

\begin{minipage}[m]{0.12\linewidth}
\centerline{\small{0.01}}
\end{minipage}
\begin{minipage}[m]{0.43\linewidth}
\centerline{~~\includegraphics[width=1.05\linewidth]{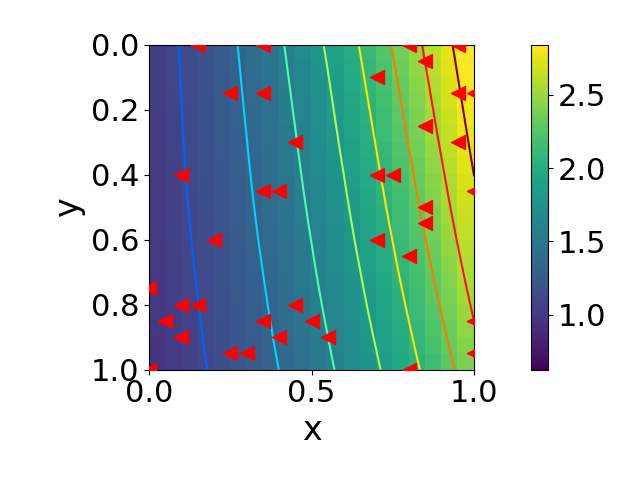}}
\end{minipage}
\begin{minipage}[m]{0.43\linewidth}
\centerline{~~\includegraphics[width=1.05\linewidth]{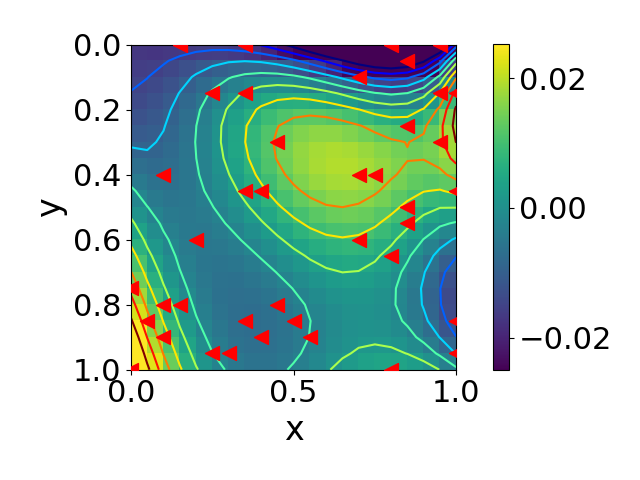}}
\end{minipage}\\

\begin{minipage}[m]{0.12\linewidth}
\centerline{\small{0.05}}
\end{minipage}
\begin{minipage}[m]{0.43\linewidth}
\centerline{~~\includegraphics[width=1.05\linewidth]{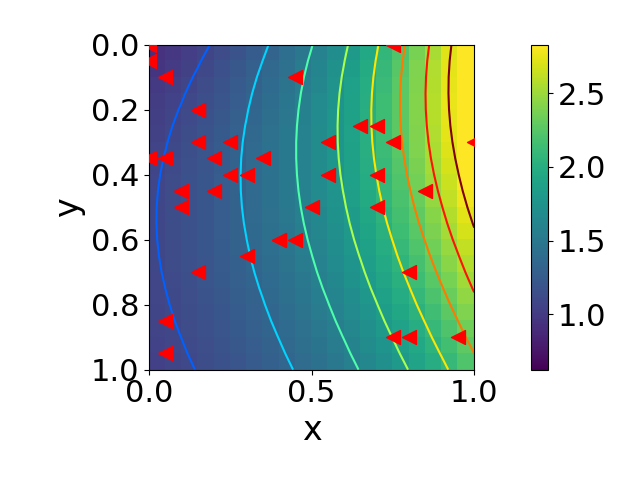}}
\end{minipage}
\begin{minipage}[m]{0.43\linewidth}
\centerline{~~\includegraphics[width=1.05\linewidth]{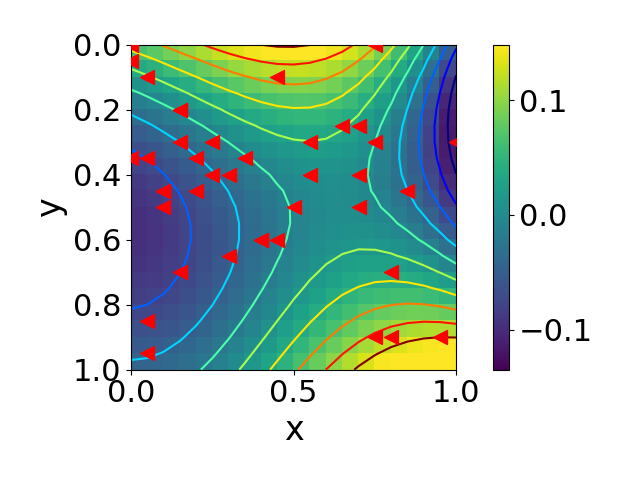}}
\end{minipage}\\
\caption{The estimations for the spatially-varying conductivity from sparse observations of noisy data and the corresponding error. The 40 randomly sampled grid points, where velocity and temperature data are observed, are labeled with red triangles.}% Each noisy observation is generated by multiplying the simulated physical quantities with a uniformly distributed scalar in $[1-\epsilon, 1+\epsilon]$. }
\label{fig:eg2_noisy}
\end{figure}

We also investigate the robustness of our method. To this end, we add a multiplicative noise sampled uniformly from $[-\epsilon, \epsilon]$ to each observation independently. The results after 100 optimization steps are shown in \Cref{fig:eg2_noisy}. We observe that the error in the conductivity estimation increases as the noise level rises, but the estimation remains very accurate given the noise level of $\epsilon=0.05$ and the sparse observations. This implies that our approach is robust to noise. 

\subsection{Learning Spatially-Varying Viscosity in Passive Transport Equations}\label{SubSec:numerical3}

% Figure

\begin{figure}[htpb] 

%\hspace{.01\linewidth}

\begin{minipage}[m]{0.32 \linewidth}
\centerline{\small{reference}~~}
\end{minipage}
\begin{minipage}[m]{0.32\linewidth}
\centerline{\small{estimation}~~~~~~}
\end{minipage}
\begin{minipage}[m]{0.32\linewidth}
\centerline{\small{difference}~~~~~~}
\end{minipage}\\

\begin{minipage}[m]{0.32\linewidth}
\centerline{\includegraphics[width=1.05\linewidth]{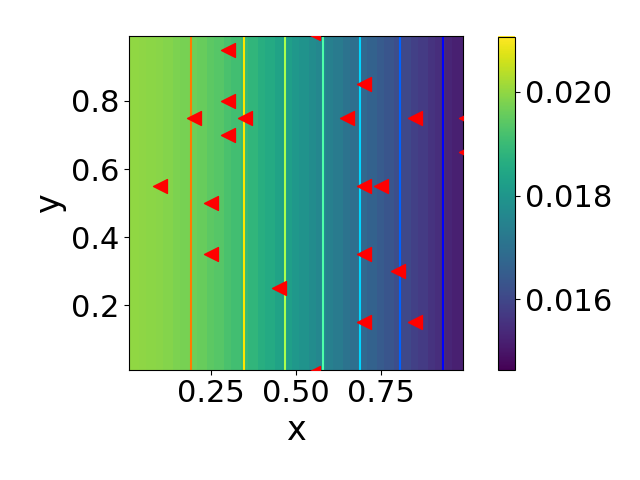}}
\end{minipage}
\begin{minipage}[m]{0.32\linewidth}
\centerline{\includegraphics[width=1.05\linewidth]{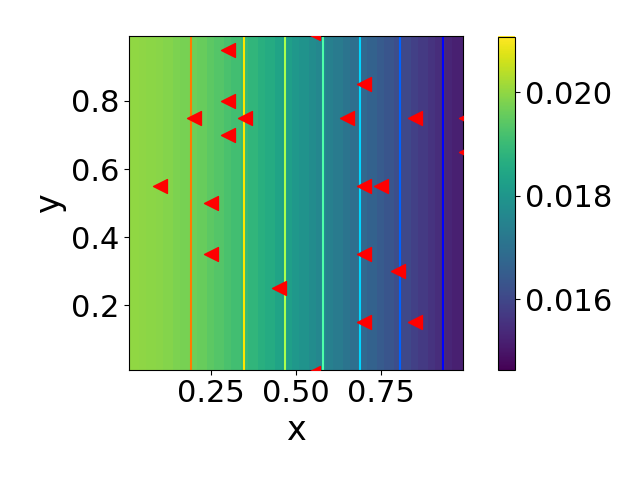}~~~~}
\end{minipage}
\begin{minipage}[m]{0.32\linewidth}
\centerline{\includegraphics[width=1.2\linewidth]{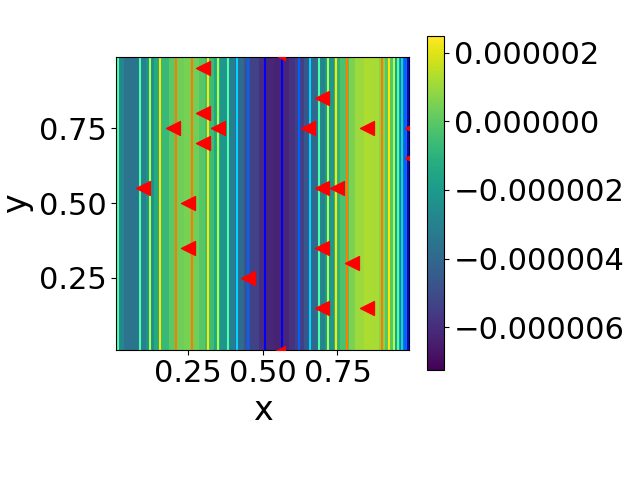}}
\end{minipage}\\
\caption{The reference viscosity, the estimated viscosity by the DNN, and the estimation error. The 22 randomly sampled grid points, where velocity data are observed, are labeled with red triangles. }\label{fig:eg3_covid}
\end{figure}

We consider an application of our method to estimate the spatially-varying viscosity from observations of passive particles. We assume that the trajectories of a passive particle are partially observed. The governing equations for the velocities of the passive particle are
\begin{align*}
    \frac{\partial w_1}{\partial t} &= \kappa_1(u - w_1) + q_1\\
    \frac{\partial w_2}{\partial t} &= \kappa_2(v - w_2) + q_2
\end{align*}
where $(u,v)$ are the velocity field from \Cref{eq:mass_and_momentum}, $(\kappa_1, \kappa_2)$ quantify the velocity-dependent accelerations of the passive particle, $(w_1, w_2)$ are the velocities of the passive particle, and $(q_1, q_2)$ are the additional body accelerations of the passive particle. We assume that $(w_1, w_2)$ are partially observed  and we want to estimate a spatially-varying viscosity field $\nu(x,y)$. This problem appears in many applications such as the modeling of nasal drug delivery \cite{basu2020numerical}, where $(u,v)$ represent the airflow velocity and $(w_1, w_2)$ represent the droplet velocity. The observations are simulated with $\rho=1$, $\kappa_1=1$, $\kappa_2=1$, and kinematic viscosity
\[
    \nu(x,y) = 0.01 + \frac{0.01}{1+x^2}
\]
In this example, we consider a layered model for $\nu(x,y)$: the estimated viscosity at a given location depends only on the $x$ coordinate. We found that the current data are not sufficient for estimating a viscosity field that depends on both $x$ and $y$ coordinates. The reference viscosity, the estimated viscosity by the layered DNN model, and the estimation error after 100 optimization steps are summarized in \Cref{fig:eg3_covid}.

\section{Discussion} \label{Sec:discussion}

Despite the generality of our approach, there are some limitations to our current work. \begin{enumerate}

\item  The memory cost is large. Due to the nature of reverse-mode automatic differentiation, we need to save all the intermediate results. This poses a big challenge when the application requires high resolution for numerical simulations. One remedy is to consider distributed computing. For example, we can use Message Passing Interface (MPI) techniques to scale the problem by utilizing multiple processors and computer nodes \cite{gropp1996high}.
This is under development for the ADCME library. 

\item Optimization with DNNs leads to a nonconvex problem, which is difficult for gradient-based optimization algorithms as local minima are inevitable. One approach is to impose some prior knowledge to the DNNs. For example, in \ref{SubSec:numerical3} we considered a layered model for the viscosity field. Although this does not solve the non-convex problem, we shrink the space of possible solutions and therefore make the inverse problem better conditioned.

\item It is difficult to determine whether an inverse problem is ill-posed before solving the inverse problem. Multiple distinct physical property  fields may produce similar or identical observations. We plan to develop diagnostic guidance for determining when an inverse problem is ill-posed in the formulation of our approach. 
    
\end{enumerate}
% limitation
% memory + computational cost large
% neural network nonconvex problem --> initial guess 
% uniqueness 

\section{Conclusion} \label{Sec:conclusion}

We have proposed a novel and general approach for solving inverse problems for the steady-state Navier-Stokes equations. In particular, we consider estimating spatially-varying physical properties (e.g., viscosity and conductivity) in coupled systems from (partially observed) state variables. The key is to express the numerical simulation using a computational graph and implement the forward computation using operators (nodes in the computational graph) that can back-propagate gradients. Then, the gradients of the loss functions with respect to the unknown parameters can be extracted automatically. We approximate the unknown physical field using a DNN and calibrate its weights and biases using a gradient-based optimization algorithm. Computing the gradients requires back-propagating gradients through both the numerical PDE solvers and the DNNs. 

Our major finding is that the DNNs provide regularization compared to pixel-wise approximations in the case of small and indirect data (i.e., partially observed state variables). We demonstrate the effectiveness and versatility of our approach with three different inverse modeling problems that involve the steady-state Navier-Stokes equations. Our implementation leverages the following two open access libraries, both of which can be easily generalized and applied to other inverse problems:
\begin{enumerate}
    \item ADCME.jl\footnote{\url{https://github.com/kailaix/ADCME.jl}}: automatic differentiation backend;
    \item AdFem.jl\footnote{\url{https://github.com/kailaix/AdFem.jl}}: a collection of numerical simulation operators.
\end{enumerate}

\section*{Acknowledgement}
This research was supported by the U.S. Department of Energy, Office of Advanced Scientific Computing Research under the Collaboratory on Mathematics and Physics-Informed Learning Machines for Multiscale and Multiphysics Problems (PhILMs) project, PhILMS grant DE-SC0019453.

This work was performed in part during an internship of T. F. at Ansys, Inc.
We acknowledge Ansys for support and Rishikesh Ranade, Haiyang He, Amir Maleki, Jan Heyse, and Wentai Zhang on the Chief Technology Officer team for helpful suggestions.

We thank the anonymous reviewers for their constructive comments on an earlier version of this paper.

\bibliographystyle{aaai}
\bibliography{ref}

% \pagebreak

%\noindent Congratulations on having a paper selected for inclusion in an AAAI Press proceedings or technical report! This document details the requirements necessary to get your accepted paper published using PDF\LaTeX{}. If you are using Microsoft Word, instructions are provided in a different document. AAAI Press does not support any other formatting software. 

%The instructions herein are provided as a general guide for experienced \LaTeX{} users. If you do not know how to use \LaTeX{}, please obtain assistance locally. AAAI cannot provide you with support and the accompanying style files are \textbf{not} guaranteed to work. If the results you obtain are not in accordance with the specifications you received, you must correct your source file to achieve the correct result. 

%These instructions are generic. Consequently, they do not include specific dates, page charges, and so forth. Please consult your specific written conference instructions for details regarding your submission. Please review the entire document for specific instructions that might apply to your particular situation. All authors must comply with the following:

\end{document}